\documentclass{article}
\usepackage{amssymb}
\usepackage{graphicx}
\usepackage{amsmath}

\newtheorem{theorem}{Theorem}

\def\bbr{\mathbb{R}}
\def\bbc{\mathbb{C}}

\begin{document}
\baselineskip=18pt
\renewcommand{\thefootnote}{}

\title{Counting Primes, Groups and Manifolds}
\author{Dorian Goldfeld$^1$, Alexander Lubotzky$^2$,\\ Nikolay Nikolov$^3$ and L\'{a}szl\'{o} Pyber$^4$}
\maketitle

\textbf{Addresses:}

\textbf{1.} Department of Mathematics, Columbia University, New York, NY 10027, U.S.A.,

\textbf{2.} Einstein Institute of Mathematics, Hebrew University, Jerusalem 91904, Israel,

\textbf{3.} Tata Institute for Fundamental Research, Colaba, Mumbai 400005, India,

\textbf{4.} A. R\'{e}nyi Institute of Mathematics, Re\'{a}ltanoda ut. 13-15, H-1053, Budapest, Hungary.

\textbf{Corresponding author:} Alexander Lubotzky, tel +972 2 6584387, e-mail: alexlub@math.huji.ac.il .

\textbf{Pages: 11 (double spaced).}

\textbf{Wordcount: abstract 150, main text 2200.}

\textbf{ PNAS Classification.} Physical sciences: Mathematics.

\emph{2000 Mathematics Subject Classification 20H05; 22E40}

\newpage

\section*{Abstract}
Let $\Lambda = \mathrm{SL}_2(\Bbb Z)$ be the modular group
and let $c_n(\Lambda)$ be the number of congruence subgroups of $\Lambda$
of index at most $n$. We prove that  $\lim\limits_{n\to \infty}
\frac{\log c_n(\Lambda)}{ (\log n)^2/\log\log n} =
\frac{3-2\sqrt{2}}{4}.$ The proof
is based on the Bombieri-Vinogradov `Riemann hypothesis on the average' and on the solution of a new type
of extremal problem in combinatorial number theory. Similar surprisingly sharp estimates are obtained for the subgroup growth of lattices in higher rank semisimple Lie groups. If $G$ is such a Lie group and $\Gamma$ is an irreducible lattice of $G$ it turns out that the subgroup growth of $\Gamma$ is independent of the lattice and depends only on the Lie type of the direct factors of $G$. It can be calculated easily from the root system. The most general case of this result relies on the Generalized Riemann Hypothesis but many special cases are unconditional. The proofs use techniques from number theory, algebraic groups, finite group theory and combinatorics.

\newpage

\section*{Statement of results: arithmetic groups}
Let $n$ be a large integer, $\Gamma$ a finitely generated group
and $M$ a Riemannian manifold.  Denote by $\pi(n)$ the number of
primes less or equal to $n$, $s_n(\Gamma)$ is the number of
subgroups of $\Gamma$ of index at most $n$ and $b_n (M)$ is the
number of covers of $M$ of degree at most $n$.  The aim of this
note is to announce results which show that in some
circumstances, these three seemingly unrelated functions are very
much connected.  This happens, for example, when $\Gamma$ is an
arithmetic group, in which case it is also the fundamental group
of a suitable locally symmetric finite volume manifold $M$.  The
studies of $s_n(\Gamma)$ and $b_n (M)$ are then almost the same.
Moreover, if $\Gamma$ has the congruence subgroup property then
estimating $s_n(\Gamma)$ boils down to counting congruence
subgroups of $\Gamma$. The latter is intimately related to the
classical problem of counting primes.
To present our results we need more notation.
 
Let $G$ be an absolutely
simple, connected, simply connected algebraic group defined over a
number field $k$. For a finite subset of valuations of $k$
including all the archimedean ones, let $\mathcal{O}_S$ denote the
ring of $S$-integers of $k$ and set $\Gamma =G(\mathcal{O}_S)$. A
subgroup $H \leq \Gamma$ is called a congruence subgroup if there
is some ideal $I \vartriangleleft \mathcal{O}_S$ such that $H$
contains the kernel of the homomorphism $\Gamma \rightarrow
G(\mathcal{O}_S/I)$.

Let $c_n(\Gamma)$ denote the number of congruence subgroups of
index at most $n$ in $\Gamma$. The counting of congruence
subgroups in arithmetic groups has already played a role in the
proof of one of the main results of the theory of subgroup growth:
A finitely generated residually finite group $\Gamma$ has
polynomial subgroup growth (i.e. $s_n(\Gamma) = n^{O(1)}$) if and
only if $\Gamma$ is virtually solvable of finite rank (cf.
\cite{LS} and the references therein).  That theorem required
only a weak lower bound on the number congruence subgroups.  In
\cite{lub} Lubotzky proved a more precise result: there exist
numbers $a,b$ depending on $G,k$ and $S$, such that$^*$
\footnotetext{$^*$ The lower bound depended on GRH at the time but
was made unconditional in \cite{GLP}} \[ n^{\frac{a \log n}{\log
\log n}} \leq c_n(\Gamma) \leq n^{\frac{b \log n}{\log \log n}},\]
and, moreover the sequence $s_n(\Gamma)$ has much faster growth
(at least $n^{\log n}$) if the congruence subgroup property fails
for $G$. Below we determine the precise rate of growth of
$c_n(\Gamma)$. (All logarithms are in base $e$.)\medskip

Let $X$ be the Dynkin diagram of the split form of $G$ (e.g. $X=A_{n-1}$ if $G=\mathrm{SU}_n$).
Let $h$ be the Coxeter number of the root system $\Phi$ corresponding to $X$ (it is the order of the Coxeter element of the Weyl group of $X$). Then $h=\frac{|\Phi|}{l}$ where $l= \mathrm{rank}_\mathbb{C}(G)= \mathrm{rank}(X)$, and for later use define $R:=h/2$. Let
\[\gamma(G)=\frac{ (\sqrt{h(h+2)}-h)^2}{4h^2}. \]

Let GRH denote the Generalized Riemann Hypothesis for Artin-Hecke
$L$-functions of number fields as stated in \cite{Weil}. The GRH
implies in particular:

Let $k$ be a Galois number field of degree $d$ over the rationals and
let $q$ be a prime such that the cyclotomic field of $q$-th roots of unity is disjoint from $k$.  Denote by $\pi_k(x, q)$ the number of primes
$p$ with $p \le x, \; p \equiv 1(\mod q)$ and $p$ splits
completely at $k$.  Then 
\[
\left| \pi_k (x, q) - \frac{x}{d \phi(q) \log x}\right| < Cx^{\frac{1}{2}} \log x \log q 
\]
for some constant $C=C(k)>0$ depending only on $k$ (a more precise bound is given in \cite{mms}). \bigskip

The lower bound for the limit in the following Theorem was proved in \cite{GLP} and the upper bound in \cite{LN}:

\begin{theorem}\label{arithmetic} Let $G$, $\Gamma$ and $\gamma(G)$ be as defined above. Assuming GRH we have
\[ \lim_{n \rightarrow \infty} \frac{\log c_n(\Gamma)}{(\log n)^2/ \log \log n}= \gamma(G),\] and moreover, this result is unconditional if $G$ is of inner type (e.g. $G$ splits) and $k$ is either an abelian extension of $\mathbb{Q}$ or a Galois extension
of degree less than 42. \end{theorem} \medskip

An interesting aspect
of this theorem is not only that the limit exists but that it is completely
independent of $k$ and $S$, and depends only on $G$. While the
independence on $S$ is a minor point and can be proved directly, the only way
we know to prove the independence on $k$ is by applying the whole machinery
of the proof.

In \cite{GLP} the crucial special case of $\Gamma=\mathrm{SL}_2(\mathcal{O}_S)$ is proved in full. There we have $\gamma(\mathrm{SL}_2)=\frac{1}{4}(3-2\sqrt{2})$. The lower bound follows using the Bombieri-Vinogradov Theorem \cite{bomb}
and the upper bound by a massive new combinatorial analysis.

\subsection*{Lattices}

Let $H$ be a connected \emph{characteristic 0} semisimple group. By this we mean that $H=\prod_{i=1}^r G_i(K_i)$ where for each $i$, $K_i$ is a local field of characteristic 0 and $G_i$ is a connected simple algebraic group over $K_i$. We assume throughout that none of the factors $G_i(K_i)$ is compact (so that
$\mathrm{rank}_{K_i}(G_i)\geq 1$). Let $\Gamma$ be an irreducible lattice of $H$, i.e. for every infinite normal subgroup $N$ of $H$  the image of $\Gamma$ in $H/N$ is dense there.

Assume now that \[\mathrm{rank}(H):=\sum_{i=1}^r \mathrm{rank}_{K_i}(G_i)\geq 2.\]
By  Margulis' Arithmeticity Theorem
(\cite{margulis}) every irreducible lattice $\Gamma$ in $H$ is arithmetic. Also the split forms of the factors $G_i$ of $H$ are necessarily of the same type and we set $\gamma(H):=\gamma(G_i)$.

Moreover, a famous conjecture of Serre (\cite{serre}) asserts that
such a group $\Gamma$ has the congruence subgroup property. It has been proved in many cases.
This enables us to prove:
\begin{theorem}\label{lat}
Assuming  GRH and Serre's conjecture, then for every non-compact
higher rank characteristic 0 semisimple group $H$ and every irreducible lattice
$\Gamma$ in $H$ the limit \[ \lim\limits_{n\to \infty} \;
\frac{\log s_n(\Gamma)}{(\log n)^2/\log\log n}\] exists and equals
$\gamma(H)$, i.e. it is independent of the lattice $\Gamma$.

Moreover the above holds unconditionally if $H$ is a simple
connected Lie group not locally isomorphic to $D_4(\mathbb{C})$
and $\Gamma$ is a non-uniform lattice in $H$ (i.e. $H/\Gamma$ is
non-compact).
\end{theorem}

Theorem 2 shows, in particular, some algebraic similarity between
different lattices $\Gamma$ in the same Lie group $G$.  This is an
addition to other results in the theory e.g. Furstenberg's theorem
showing that the boundaries of all such $\Gamma$'s are the same or
Margulis super-rigidity, which shows that the finite dimensional
representation theory of the different $\Gamma$'s in the same $G$
are similar.  (cf \cite{margulis} and the references therein).

We point out the following geometric reformulation of the special case:
\begin{theorem}\label{t'} Let $H$ be a simple connected Lie group of $\mathbb{R}$-rank $\geq 2$ which is not locally isomorphic to $D_4(\mathbb{C})$. Put $X=H/K$ where $K$ is a maximal compact subgroup of $H$. Let $M$ be a finite volume non-compact manifold
covered by $X$ and let $b_n (M)$ be the number of covers of $M$ of
degree at most $n$.  Then $\lim\limits_{n\to\infty} \;
\frac{\log b_n(M)}{(\log n)^2/\log\log n}$ exists, equals $\gamma(H)$
and is independent of $M$.
\end{theorem}

It is interesting to compare Theorems 2 and 3 with the results of
Liebeck-Shalev\cite{ls}, and M\"{u}ller-Puchta \cite{MP}: If
$H=\mathrm{SL}_2(\mathbb{R})$ and $\Gamma$ is a lattice in $H$
then $\lim\limits_{n\to \infty} \;\frac{\log s_n(\Gamma)}{\log n!}=-\chi(\Gamma)$, where
$\chi$ is the Euler characteristic. \bigskip

We finally mention a conjecture and a question:  Let $X$ be the
symmetric space associated with a simple Lie group $H$ as in
Theorem 3.  Denote by $m_n(X)$ the number of manifolds covered by
$X$ of volume at most $n$.  By a well known result of Wang \cite{Wang},
this number is finite unless $H$ is locally isomorphic to
$SL_2(\bbr)$ or $SL_2(\bbc)$.

\bigskip
\noindent {\bf Conjecture.} \  If $\mathbb{R}\textrm{-rank}(H) \ge 2$ then
\[
\lim \limits_{n\to \infty} \; \frac{\log m_n(X)}{(\log n)^2/\log\log n} = \gamma(H).
\]

\noindent {\bf Question:} Estimate $m_n(X)$ for the case of $H$
having $\bbr$-rank equal to one.  For $H=\mathrm{SO}(n, 1)$ the results of
\cite{BGLM} suggest that $\lim \limits_{n\to \infty} \; \frac{\log m_n
(H) }{\log n!}$ may exist, but we do not have any clue what it
could be.
\section*{Proofs: the lower bound}
We shall illustrate the main idea of the proof with $\Gamma= \mathrm{SL}_d(\mathbb{Z})$ and refer to \cite{GLP} for the full details.

Choose any $\rho \in (0,\frac{1}{2})$. For $x>>0$ and a prime $q<x$ let $P(x,q)$ be the set of primes $p\leq x$ such that $p \equiv 1$ mod $q$. Let $L(x,q)=|P(x,q)|$ and $M(x,q)= \sum_{p\in P(x,q)} \log p$. Then the Bombieri-Vinogradov Theorem \cite{bomb} ensures the existence of a prime $q \in (\frac{x^\rho}{\log x},x^\rho)$ such that
\[ L(x,q) = \frac{x}{\phi(q) \log x}+O\left(\frac{x}{\phi(q) (\log x)^2}\right) ;\quad M(x,q) = \frac{x}{\phi(q)}+O\left(\frac{x}{\phi(q) (\log x)^2}\right).\]
Put $L:=L(x,q)$ and $M:=M(x,q)$.

By strong approximation  (cf. \cite{LS}, Window 9) $\Gamma$ maps onto \\
$G_P:=\prod_{p\in P(x,q)}\mathrm{SL}_d(\mathbb{F}_p)$. Let $B(p)$
be the subgroup of upper triangular matrices of
$\mathrm{SL}_d(\mathbb{F}_p)$ and set
\[B_P:=\prod_{p\in P(x,q)}B(p).\] The group $B_P$ maps onto the
diagonal $\prod_p (\mathbb{F}^*_p)^{d-1}$ which in turn maps onto
$\mathbb{F}_q^{(d-1)L}$. For fixed $\sigma \in (0,1)\cap
\frac{1}{L(d-1)}\mathbb{N}$ the latter vector space has about
$q^{\sigma (1-\sigma)(d-1)^2 L^2}$ subgroups of index $q^{\sigma
(d-1)L}$ (see Proposition 1.5.2 in \cite{LS}), each giving rise to
a subgroup of index $n=[G_P:B_P]q^{\sigma (d-1)L}$ in $\Gamma$.
Now \\ $\log [G_P:B_P]\sim d(d-1)M/2$ as $x \rightarrow \infty$
and after some algebraic manipulations we obtain that for this
chosen value of $n$
\[\frac{\log c_n(\Gamma)}{(\log n)^2/ \log \log n} \geq\frac{\sigma (1-\sigma)\rho (1-\rho)}{\left(\sigma \rho +R\right)^2} -o(1), \quad (x \rightarrow \infty) \] where in our case $R=d/2$. As shown in \cite{GLP} \S 3 the maximum value of the above expression for $\sigma, \rho \in (0,1)$ is precisely $\gamma (G)=\frac{ (\sqrt{R(R+1)}-R)^2}{4R^2}$ and is achieved for $\sigma_0=\rho_0= \sqrt{R(R+1)}-R$. By taking $x$ sufficiently large we can choose $\sigma \in (0,1)\cap \frac{1}{L(d-1)}\mathbb{N}$ to be arbitrarily close to $\sigma_0$, and take $\rho=\rho_0$. This proves the lower bound.  \medskip

The reason for invoking the GRH in Theorem \ref{arithmetic} is that in the general case we need an equivalent of the Bombieri-Vinogradov theorem for $k$ in place of $\mathbb{Q}$. The work of M.R. Murty and V.K. Murty \cite{murty} gives an analogue of it for number fields but their result is weaker in general. It suffices for our needs when, for example $k/\mathbb{Q}$, is an abelian extension.
\subsection*{The upper bound}
The proof of the upper bound in \cite{LN} is inspired by the special case solved in \cite{GLP} and has two parts:

I. A reduction to an extremal problem for abelian groups, and

II. Solving this extremal problem (Theorem \ref{ab} below).

\textbf{Part I:}

The subgroup structure of the groups $\mathrm{SL}_2(\mathbb{F}_p)$ is completely known. Using
this it is shown in \cite{GLP} that Theorem \ref{arithmetic} for $\mathrm{SL}_2(\mathbb{Z})$ is equivalent to the
following extremal result on counting subgroups of abelian groups:

Let $C_m$ denote the cyclic group of order $m$.
For all pairs $\mathcal{P}_-$ and $\mathcal{P}_+$ of disjoint sets of primes, let \[ f(n):= \max \left \{ s_r(X)\ \left \vert  \quad X= \prod_{p \in \mathcal{P}_-} C_{p-1} \times  \prod_{p \in \mathcal{P}_+} C_{p+1} \right. \right \}, \]
where the maximum is taken over all sets $\mathcal{P}_-, \mathcal{P}_+$ and $r \in \mathbb{N}$ such that \\ $n \geq r \prod_{p \in \mathcal{P}} p$, (here $\mathcal{P}=\mathcal{P}_- \cup \mathcal{P}_+$).
\begin{theorem} We have \[
\limsup_{n \rightarrow \infty} \frac{ \log c_n(\mathrm{SL}_2(\mathbb{Z}))}{(\log n)^2 / \log \log n} = \limsup_{n \rightarrow \infty} \frac{ \log f(n)}{(\log n)^2 / \log \log n}.\]
\end{theorem}

By contrast there is no such precise description of the subgroup structure even for $\mathrm{SL}_n(\mathbb{F}_p)$. Still, surprisingly, the proof of the general upper
bound reduces to a similar extremal problem for abelian groups using some
ideas of \cite{GLP}, \cite{LP} and the following Theorem which is the
main new ingredient in \cite{LN}.

Let $X(\mathbb{F}_q)$ be a finite quasisimple group of Lie type
$X$ over the finite field $\mathbb{F}_q$ of characteristic $p>3$.
For a subgroup $H$ of $X(\mathbb{F}_q)$ define
\[ t(H)=\frac{ \log[X(\mathbb{F}_q):H]}{\log|H^{\diamondsuit}|},\]
 where $H^{\diamondsuit}$ denotes the maximal abelian
 quotient of $H$ whose order is coprime to $p$. Set $t(H)=\infty$ if $|H^{\diamondsuit}|=1$.
\medskip

Recall that $R=R(X)=h/2$ where $h$ is the Coxeter number of the root system of the \textbf{split} Lie type corresponding to $X$.

\begin{theorem}\label{t4} Given the Lie type $X$ then
\[ \liminf_{q \rightarrow \infty} \ \min \left \{ t(H)\ |\quad H\leq X(\mathbb{F}_q)\ \right \} \geq R.\]
\end{theorem}

The proof of this theorem does not depend on the classification of the finite simple groups, we use
instead the work of Larsen and Pink \cite{lpink} (which is a classification-free version of a result of Weisfeiler \cite{We}), and Liebeck, Saxl
and Seitz \cite{lss} (the latter for groups of exceptional type).

\textbf{Part II:}

Once Part I is proved, the argument reduces to an extremal problem on abelian groups:
\begin{theorem}\label{ab} Let $d$ and $R$ be fixed positive numbers. Suppose $A=C_{x_1}\times C_{x_2}\times \cdots \times C_{x_t}$ is an abelian group such that the orders $x_1,x_2,...,x_t$ of its cyclic factors do not repeat more than $d$ times each. Suppose that $r|A|^R\leq n$ for some positive integers $r$ and $n$. Then as $n,r$ tend to infinity we have
\[s_r(A)\leq n^{(\gamma +o(1))\frac{\log n}{\log \log n}},\]
where $\gamma=\frac{ (\sqrt{R(R+1)}-R)^2}{4R^2}$.
\end{theorem}

The starting point of the proof of this theorem in \cite{GLP} is a well-known formula for counting subgroups of finite abelian groups (see \cite{Bu}). We refer the reader to \cite{GLP} for the details which are too complicated to be given here.

\subsection*{Acknowledgments} In their work the authors were supported by
grants from the National Science Foundation (Goldfeld and Lubotzky), the Israel Science Foundation and the US-Israel Binational Science Foundation (Lubotzky) and Hungarian National Foundation for Scientific Research, Grant T037846 (Pyber). While this research was carried out Nikolov held a
Golda-Meir Postdoctoral Fellowship at the Hebrew University of
Jerusalem.

\end{document}